\documentclass[A4paper,11pt]{article}
\usepackage{latexsym}
\usepackage{mathrsfs}
\usepackage{amssymb}
\usepackage{amscd}
\usepackage[dvips]{graphicx} 
\usepackage{xcolor}
\usepackage{color}

\newtheorem{theorem}{Theorem}

\newtheorem{proposition}{Proposition}
\newenvironment{proof}[1][Proof]{\textbf{#1.} }{\ \rule{0.5em}{0.5em}}
\date{}
\long\def\symbolfootnote[#1]#2{\begingroup%
	\def\thefootnote{$\;$}\footnote[#1]{$^*$#2}\endgroup}
\begin{document}
	
	\title{Partitions and point-finite covers \\of Baire spaces}
	\author{Joanna Jureczko}
\maketitle

\symbolfootnote[2]{Mathematics Subject Classification: Primary 03C25, 03E35, 03E55, 54E52}.

	\hspace{0.2cm}

\begin{abstract}
	In this paper we show that the existence of Kuratowski partitions of Hausdorff Baire space is equivalent to the existence of point-finite covers of the same space.
	
\end{abstract}

\section{Introduction}

The existence of Kuratowski partitions is related to the old question posed by K. Kuratowski, in 1935, see \cite{KK},  (following the results of Lusin published in 1912, (\cite{NL})),  who posed the problem  whether a function $f \colon X \to~Y$, (where $X$ is completely metrizable and $Y$ is metrizable), such that each preimage of an open set of $Y$ has the Baire property, is continuous apart from a meager set.

The problem has been intensively studied in the 70's of the last century among others by Solovay (unpublished results) and Bukovsky, (\cite{LB}).
Further historical information and results concerning Kuratowski partitions can be found in \cite{JJ1}.

In 1979,  in the same issue of Bulletin of Polish Academy of Sciences (Bull. Ac. Pol.: Math.  27 (1979)) there were published three papers concerning this problem.
 In the first two \cite{EFK, EFK1} A. Emeryk, R. Frankiewicz, and W. Kulpa  demonstrated that Kuratowski's problem is equivalent to asserting the existence of partitions of completely metrizable spaces into meager sets such that the union of each subfamily  of this partition has the Baire property, (such a partition is called a \textit{Kuratowski partition}) and showed that if $\mathcal{F}$ is a partition of \v Cech complete space of weight $\leqslant 2^\omega$, then there exists a family $\mathcal{F}' \subset \mathcal{F}$ such that $\bigcup \mathcal{F}'$ has not the Baire property. 
 In the third paper, \cite{BCGR}, J. Brzuchowski, J. Cicho\'n, E. Grzegorek and C. Ryll-Nardzewski proved the similar result for point-finite covers instead of partitions.

 The main result of this paper is to show that the existence of point-finite covers of a Hausdorff Baire space into meager sets being  completely additive with respect to the Baire property is equivalent to the existence of a partition of this space into meager sets with the same additive property, (Theorem 1 and Theorem 2).
 It is worth adding that it is impossible to obtain the result similar to Theorem 2 for the whole space, because, as shown in \cite{JJ} there exists a non-complete metric space with Kuratowski partition for which its completion has no Kuratowski partition.
 
 Having such a result we can give the adequate results for point-finite covers which were previously given for partitions, (see \cite{JJ, JJ3, JJ4}). In this paper we show that the result concerning the existence of measurable cardinals (presented in \cite{FK} for partitions) is true for point-finite covers.  Other examples are given in \cite{FJW} where it is shown that in the Marczewski and Laver structure and structure with Ellentuck topology there are no Kuratowski partitions. Recently, it was shown that the similar results can be obtained in these structures for point-finite covers instead of partitions, (see \cite{JJ2}). 
 
 For definitions and  facts not cited here we refer to e.g. \cite{RE, KK1} (topology) and \cite{TJ} (set theory).

\section{Definitions and previous results}
Throughout the paper, we assume that $X$ is a Hausdorff Baire space, i.e. a Hausdorff space in which the Baire theorem holds.
\\\\
   A set $U \subseteq X$ has \textit{the Baire property} iff there exists an open set $V \subset X$ and a meager set $M \subset X$ such that $U = V \triangle M$, where $\triangle$ represents the symmetric difference of sets.
   \\
   \\
   A set $A \subseteq X$ is \textit{co-meager} in $X$ iff $X \setminus A$ is meager in $X$.
   \\\\
   \textbf{Fact 1. (folklore)} Let $X$ be a  Baire space. Let $A \subseteq X$ be a co-meager set in $X$ and let $C \subseteq X$ be a meager set in $X$. Then the set $A \cap C$ is meager in $A$.    
\\
\\
A family $\mathcal{A} \subset P(X)$ is \textit{completely additive with respect to property $\mathcal{P}$}, shortly $\mathcal{A}$ is $CA(\mathcal{P})$, iff $\bigcup \mathcal{A}'$ has  property $\mathcal{P}$ for all $\mathcal{A}' \subseteq \mathcal{A}$.
\\\\ Let $X$ be a fixed space. A partition $\mathcal{F}$ of $X$ into meager subsets of $X$ is called a \textit{Kuratowski partition} iff $\mathcal{F}$ is $CA(Baire)$. We define and assume in the whole paper that $$\kappa = \min\{|\mathcal{F}|\colon \mathcal{F} \textrm{ is Kuratowski partition of }X \}.$$
For a given cardinal $\kappa$ we enumerate
$\mathcal{F} = \{F_\alpha \colon \alpha < \kappa\}$. In the whole paper we  can assume that $\kappa$ is a regular cardinal.  This is the standard assumption,(see e.g. \cite{TJ}). If $\kappa$ was singular then $cf(\kappa)$ would be the minimal one. Indeed, each union of less than $\kappa$ many sets is meager, which follows from Union Theorem \cite[p.82]{KK1}. If it would be non-meager and hence contain an open set, would give Kuratowski partition of size less than $\kappa$. By Baire Theorem, $\kappa$ is uncountable.
\\
\\
With any Kuratowski partition
$\mathcal{F} = \{F_\alpha \colon \alpha < \kappa\}$ indexed by  $\kappa$,  we may associate an ideal 
$$I_\mathcal{F} = \{A \subset \kappa \colon \bigcup_{\alpha \in A} F_\alpha \textrm{ is meager}\},$$
which we call  \textit{$K$-ideal}. Obviously, $I_\mathcal{F}$ is non-principal and $[\kappa]^{< \kappa}~\subseteq I_{\mathcal{F}}$. 
 ${K}$-ideal is $\kappa$-complete which is the consequence of Union Theorem \cite[p.82]{KK1}. Moreover, by \cite[Fact 13]{JJ} if $\pi w(X)\leqslant 2^\omega$ then there are no Kuratowski partitions on $X$. Hence $\kappa$ is at least $\omega_1$.
\\
\\
A family $\mathcal{F}$ is called \textit{point-finite} if every point of $X$ lies in only finitely many members of $\mathcal{F}$.

\section{Main results}
\begin{theorem}
	Let $X$ be a Hausdorff space and let $A\subset X$ be a Baire subspace of $X$ which admits Kuratowski partition. Then there is a $CA(Baire)$ point-finite cover $\mathcal{G}$ of $A$ by meager sets. 
	\end{theorem}

\begin{proof}
Let $\mathcal{F}$ is a Kuratowski partition of $A$. Then $$|\{F \in \mathcal{F} \colon x \in F\}|=1$$ for any $x \in A$. Thus, $\mathcal{F}$ is a point finite cover of $A$.
\end{proof}

\begin{proposition} Suppose that $n < \omega$ and $\{Z_\alpha \colon \alpha < \kappa\}$ is a cover of a set $X$ such that $|\{\alpha < \kappa \colon x \in Z_\alpha\}| = n$ for each $x \in X$. Then there are disjoint sets $\{Y_\alpha \colon \alpha < \kappa\}$ such that
	$$Y_\alpha \subseteq Z_\alpha \textrm{ for $\alpha < \kappa$ and  } \bigcup_{\alpha< \kappa} Y_\alpha = X.$$ 
	\end{proposition}

\begin{proof}
	Let $\{Z_\alpha \colon \alpha < \kappa\}$ be a cover of a set $X$. By assumptions $$|\{\alpha < \kappa \colon x \in Z_\alpha\}| = n \textrm{ for each } x \in X.$$ For each $x \in X$ denote 
	$$A^n_x = \{\alpha < \kappa \colon x \in Z_\alpha\}.$$
	Obviously, $|A^n_x| = n$. 
	For each such $\alpha < \kappa$ we will construct a decreasing sequence of sets $\{Z^k_{\alpha} \colon 1 \leqslant k \leqslant n\}$. 
	We start with $Z^n_{\alpha} = Z_{\alpha}$. Then take 
	$$Z^{n-1}_{\alpha} = Z^n_{\alpha} \setminus \{x \in X \colon \min A^n_x = \alpha\}.$$
	Now, define $A^{n-1}_x = \{\alpha < \kappa \colon x \in Z^{n-1}_\alpha\}$. Then $|A^{n-1}_{x}| = n-1$.
	Next, take 	$$Z^{n-2}_{\alpha} = Z^{n-1}_{\alpha} \setminus \{x \in X \colon \min A^{n-1}_x = \alpha\}$$
	and define $A^{n-2}_x = \{\alpha < \kappa \colon x \in Z^{n-2}_\alpha\}$. Then $|A^{n-2}_{x}| = n-1$.
	Continuing the construction we finally define $Z^1_{\alpha}$ such that
	$|A^1_x| = |\{\alpha < \kappa \colon x \in Z^1_{\alpha}\}| =1$. 
	Thus, for each $\alpha < \kappa$ we obtain a pairwise disjoint subsets $\{Z^1_{\alpha} \colon \alpha< \kappa\}$ of $X$. Putting  $Y_\alpha = Z_\alpha^1$ we obtain our claim.  
\end{proof}

	\begin{theorem}
	Let $X$ be a Hausdorff Baire space which allows a $CA(Baire)$ point-finite cover by meager sets. Then there exists a co-meager subset $B \subseteq X$ which admits Kuratowski partition. 
	\end{theorem}
	
	\begin{proof} 
		Let $\mathcal{G}$ denote a $CA(Baire)$ point-finite cover of $X$ by meager sets.	
 	For any $n \in \omega$, consider sets
 	$$W(n) = \{x \in A \colon |\{G \in \mathcal{G} \colon x \in G\}| =n\}.$$
 	Notice that $W(n) \cap W(m) = \emptyset$, whenever $n \not = m$, for any $n, m < \omega$.
 	By Baire Theorem some $W(n)$ are non-meager. Put
 	$$K = \{n \in \omega \colon W(n) \textrm{ is non-meager}\}.$$
 	Let
 	$$B= \bigcup\{W(n) \colon n \in K\}.$$
 	Obviously, $B$ is co-meager.
 	Let 
 	$$\mathcal{G}_B = \{G \cap B \colon G \in \mathcal{G}\}$$ be a point-finite cover of $B$.
 	By Fact 1, we notice that $G\cap B \in \mathcal{G}_B$ is meager.
 	
 	Now, for each $n \in K$ and each $x \in W(n)$ define
 	 $$C^n_x = \{G \cap B \in \mathcal{G}_B \colon x \in G\}.$$
 	 Obviously, $|C^n_x| =n$.
 	 For each such $C^n_x$ we use Proposiion 1 obtaining a family of disjoint meager sets 
 	 $$\mathcal{F}_B = \{F \subseteq G \cap B \colon G \cap B \in C^n_x, n \in K, x \in W(n)\}.$$
 	 By Union Theorem \cite[p.82]{KK1} the family $\mathcal{F}_B$ is a required Kuratowski partition.  
	\end{proof}
	
	\section{Consequences}
	
	\begin{theorem}
		IF $ZFC+$ there is a Baire Hausdorff space $X$ and a $CA(Baire)$ point-cover of $X$ consisting of meager sets  is consistent, then $ZFC+$ there is a measurable cardinal is consistent as well.
			\end{theorem}
			
			\begin{proof}
Let $\kappa$ be the smallest regular cardinal of $CA(Baire)$ point-finite cover ofsome Baire space $X$ into meager sets. Choose such a cover and denote by $\mathcal{G}$.

By Theorem 2, there exists a co-meager subset $B \subseteq X$ which admits Kuratowski partition.

Let $\mathcal{F} = \{F_\alpha \colon \alpha < \kappa\}$ be a Kuratowski partition of $B$.  Consider a $K-$ideal $I_\mathcal{F}$ associated with $\mathcal{F}$. Such an ideal is $\kappa$-complete and non-principal.
			
Now, we will modify the proof presented in \cite{FK}. Define the family $P(V)$  of functions, where $V$ is the universe, i.e. $f \in P(V)$ iff $f \colon X  \to V$ and there exists an open cover  $\mathcal{U}_{f}$ of  $X$, dense in $X$ such that $f$ is constant on $U \cap F_\alpha$ for any $U\in \mathcal{U}_{f}$ and $F_\alpha \in \mathcal{F}$.
			
For any $f, g \in P(V)$ the sets $$\{x \in \bigcup \mathcal{F} \colon f(x) = g(x)\} \textrm{ and } \{x \in \bigcup \mathcal{F} \colon f(x) \in  g(x)\}$$ have the Baire property. To see this take
$W \in \{U \cap V \colon U \in \mathcal{U}_{f}, V \in \mathcal{U}_{g}\}$. Then
$$\{x \in \bigcup \mathcal{F}\{\cap W \colon f(x) = g(x)\} = \bigcup \{F_\alpha \cap W \colon \alpha \in A\}$$
and
$$\{x \in \bigcup \mathcal{F}\{\cap W \colon f(x) \in g(x)\} = \bigcup \{F_\alpha \cap W \colon \alpha \in B\}$$
for some $A, B \subseteq \kappa$.
Since $\mathcal{F}$ is Kuratowski partition of $B$ then $$\mathcal{F} \cap W = \{F_\alpha \cap W \colon F_\alpha \in \mathcal{F}\}$$ is Kuratowski partition of $W$. By Banach Localisation Theorem, both $$\{x \in \bigcup \mathcal{F} \colon f(x) = g(x)\} \textrm{ and } \{x \in \bigcup \mathcal{F} \colon f(x) \in g(x)\}$$ have the Baire property.
			
 Let $\mathbb{B}$ be the Boolean algebra of regular open subsets of $X$ and $G$-generic ultrafilter over $\mathbb{B}$. 
 
 Consider  a limit ultrapower $P(V)/G$ (in Keisler sense, see e.g. \cite{CK}) which is a model of $ZFC$. By \L{}o\'s Theorem, $\cite[p. 159]{TJ}$, we have
  $$P(V)/G \models \varphi ([f_1], [f_2], ..., [f_n])$$ iff  $$\{x \in \bigcup \mathcal{F} \colon \models \varphi (f_1(x), f_2(x), ..., f_n (x))\} \in G,$$
   where $f_1, f_2, ..., f_n \in P(V)$. 
   
   As in \cite[p. 284-288]{TJ} we can define a natural embedding $j_G \colon V \to P(V)/G$, $j_G(x) = [c_x]$, where $c_x \colon \bigcup\mathcal{F} \to V$ such that $c_x(t) = x$ for all $t \in \bigcup\mathcal{F}$.
   $P(V)/G$ is well-founded, (see \cite[p. 286]{TJ}), i.e. $$P(V)/G \models f_1\ni f_2\ni ... \ni f_n \ni ...,$$ then  $$\{x \in \bigcup\mathcal{F} \colon f_n(x)\ni f_{n+1}(x)\}$$ is comeager in $\bigcup \mathcal{F}$. 
   
   By Baire Category Theorem, there exists $x_0 \in \bigcup \mathcal{F}$
   $$f_0(x_0) \ni f_1(x_0) \ni ... \ni f_n(x_0) \ni ...$$
   As shown in \cite[p. 287]{TJ} $j_G(\kappa)> \kappa$. By \cite[p. 287]{TJ} there exists a measurable cardinal.	
			\end{proof}
	\\
	\\
 \textbf{Acknowledgments} The author is very grateful to the reviewers for their thorough reading of the text and valuable comments that allowed to avoid errors and shortcomings in the presented results.

	{\sc Joanna Jureczko}
	\\
	Wroc\l{}aw University of Science and Technology, Wroc\l{}aw, Poland
	\\
	{\sl e-mail: joanna.jureczko@pwr.edu.pl}

\end{document}